%
%

\catcode`@=11
%
---

\def\b@lank{ }

\newif\if@simboli
\newif\if@riferimenti
\newif\if@bozze

\def\bozze{\@bozzetrue\font\tt@bozze=cmtt8}
\def\og@gi{\number\day\space\ifcase\month\or
   gennaio\or febbraio\or marzo\or aprile\or maggio\or giugno\or
   luglio\or agosto\or settembre\or ottobre\or novembre\or dicembre\fi
   \space\number\year}
\newcount\min@uti
\newcount\or@a
\newcount\ausil@iario

\min@uti=\number\time
\or@a=\number\time
\divide\or@a by 60
\ausil@iario=-\number\or@a
\multiply\ausil@iario by 60
\advance\min@uti by \number\ausil@iario
\def\ora@esecuzione{\the\or@a:\the\min@uti}  
\def\makefootline{\baselineskip=24pt\line{\the\footline}
    \if@bozze\vskip-10pt\tt@bozze
             \noindent \jobname\hfill\og@gi, ore \ora@esecuzione\fi}

\newwrite\file@simboli
\def\simboli{
    \immediate\write16{ !!! Genera il file \jobname.SMB }
    \@simbolitrue\immediate\openout\file@simboli=\jobname.smb}

\newwrite\file@ausiliario
\def\riferimentifuturi{
    \immediate\write16{ !!! Genera il file \jobname.AUX }
    \@riferimentitrue\openin1 \jobname.aux
    \ifeof1\relax\else\closein1\relax\input\jobname.aux\fi
    \immediate\openout\file@ausiliario=\jobname.aux}

\newcount\eq@num\global\eq@num=0
\newcount\sect@num\global\sect@num=0
\newcount\teo@num\global\teo@num=0
\newif\if@ndoppia
\def\numerazionedoppia{\@ndoppiatrue\gdef\la@sezionecorrente{\the\sect@num}}

\def\se@indefinito#1{\expandafter\ifx\csname#1\endcsname\relax}
\def\spo@glia#1>{} 

\newif\if@primasezione
\@primasezionetrue

\def\s@ection#1\par{\immediate
    \write16{#1}\if@primasezione\global\@primasezionefalse\else\goodbreak
    \vskip\spaziosoprasez\fi\noindent
    {\sezfont #1}\nobreak\vskip\spaziosottosez\nobreak\noindent}
%

\font\sezfont=cmbx10

\def\sezpreset#1{\global\sect@num=#1
    \immediate\write16{ !!! sez-preset = #1 }   }

\def\spaziosoprasez{50pt plus 60pt}
\def\spaziosottosez{15pt}
\def\spaziotitsez{5truemm}

\def\sref#1{\se@indefinito{@s@#1}\immediate\write16{ ??? \string\sref{#1}
    non definita !!!}
    \expandafter\xdef\csname @s@#1\endcsname{??}\fi\csname @s@#1\endcsname}

\def\autosez#1#2\par{
    \global\advance\sect@num by
    1\if@ndoppia\global\eq@num=0\teo@num=0\fi
    \xdef\la@sezionecorrente{\the\sect@num}
    \def\usa@getta{1}\se@indefinito{@s@#1}\def\usa@getta{2}\fi
    \expandafter\ifx\csname @s@#1\endcsname\la@sezionecorrente\def
    \usa@getta{2}\fi
    \ifodd\usa@getta\immediate\write16
      { ??? possibili riferimenti errati a \string\sref{#1} !!!}\fi
    \expandafter\xdef\csname @s@#1\endcsname{\la@sezionecorrente}
    \immediate\write16{\la@sezionecorrente. #2}
    \if@simboli
      \immediate\write\file@simboli{ }\immediate\write\file@simboli{ }
      \immediate\write\file@simboli{  Sezione
                                  \la@sezionecorrente :   sref.   #1}
      \immediate\write\file@simboli{ } \fi
    \if@riferimenti
      \immediate\write\file@ausiliario{\string\expandafter\string\edef
      \string\csname\b@lank @s@#1\string\endcsname{\la@sezionecorrente}}\fi
    \goodbreak\vskip\spaziosoprasez
    \noindent\if@bozze\llap{\tt@bozze#1\ }\fi
      {\sezfont\the\sect@num.\hskip\spaziotitsez #2}\par\nobreak
    \vskip\spaziosottosez\nobreak\noindent}

\def\semiautosez#1#2\par{
    \gdef\la@sezionecorrente{#1}\if@ndoppia\global\eq@num=0\teo@num=0\fi
    \if@simboli
      \immediate\write\file@simboli{ }\immediate\write\file@simboli{ }
      \immediate\write\file@simboli{  Sezione ** : sref.
          \expandafter\spo@glia\meaning\la@sezionecorrente}
      \immediate\write\file@simboli{ }\fi
    \s@ection#2\par}


\def\eqpreset#1{\global\eq@num=#1
     \immediate\write16{ !!! eq-preset = #1 }     }
\def\teopreset#1{\global\teo@num=#1
     \immediate\write16{ !!! teo-preset = #1 }     }

\def\eqref#1{\se@indefinito{@eq@#1}
    \immediate\write16{ ??? \string\eqref{#1} non definita !!!}
    \expandafter\xdef\csname @eq@#1\endcsname{??}
    \fi\csname @eq@#1\endcsname}

\def\teoref#1{\se@indefinito{@teo@#1}
    \immediate\write16{ ??? \string\teoref{#1} non definita !!!}
    \expandafter\xdef\csname @teo@#1\endcsname{??}
    \fi\csname @teo@#1\endcsname}

\def\eqlabel#1{\global\advance\eq@num by 1
    \if@ndoppia\xdef\il@numero{\la@sezionecorrente.\the\eq@num}
       \else\xdef\il@numero{\the\eq@num}\fi
    \def\usa@getta{1}\se@indefinito{@eq@#1}\def\usa@getta{2}\fi
    \expandafter\ifx\csname @eq@#1\endcsname
    \il@numero\def\usa@getta{2}\fi
    \ifodd\usa@getta\immediate\write16
       { ??? possibili riferimenti errati a \string\eqref{#1} !!!}\fi
    \expandafter\xdef\csname @eq@#1\endcsname{\il@numero}
    \if@ndoppia
       \def\usa@getta{\expandafter\spo@glia\meaning
       \la@sezionecorrente.\the\eq@num}
       \else\def\usa@getta{\the\eq@num}\fi
    \if@simboli
       \immediate\write\file@simboli{  Equazione
            \usa@getta :  eqref.   #1}\fi
    \if@riferimenti
       \immediate\write\file@ausiliario{\string\expandafter\string\edef
       \string\csname\b@lank @eq@#1\string\endcsname{\usa@getta}}\fi}

\def\teolabel#1{\global\advance\teo@num by 1
    \if@ndoppia\xdef\il@numero{\la@sezionecorrente.\the\teo@num}
       \else\xdef\il@numero{\the\teo@num}\fi
    \def\usa@getta{1}\se@indefinito{@teo@#1}\def\usa@getta{2}\fi
    \expandafter\ifx\csname @teo@#1\endcsname\il@numero\def\usa@getta{2}\fi
    \ifodd\usa@getta\immediate\write16
       { ??? possibili riferimenti errati a \string\teoref{#1} !!!}\fi
    \expandafter\xdef\csname @teo@#1\endcsname{\il@numero}
    \if@ndoppia
       \def\usa@getta{\expandafter\spo@glia\meaning
       \la@sezionecorrente.\the\teo@num}
       \else\def\usa@getta{\the\teo@num}\fi
    \if@simboli
       \immediate\write\file@simboli{  Equazione
            \usa@getta :  teoref.   #1}\fi
    \if@riferimenti
       \immediate\write\file@ausiliario{\string\expandafter\string\edef
       \string\csname\b@lank @teo@#1\string\endcsname{\usa@getta}}\fi}

\def\autoreqno#1{\eqlabel{#1}\eqno(\csname @eq@#1\endcsname)
       \if@bozze\rlap{\tt@bozze\ #1}\fi}
\def\autoteono#1{\teolabel{#1}\csname @teo@#1\endcsname
       \if@bozze\rlap{\raise2ex\hbox{\tt@bozze\ #1}}\fi}

\def\autoleqno#1{\eqlabel{#1}\leqno\if@bozze\llap{\tt@bozze#1\ }
       \fi(\csname @eq@#1\endcsname)}
\def\eqrefp#1{(\eqref{#1})}
\def\numeriadestra{\let\autoeqno=\autoreqno}
\def\numeriasinistra{\let\autoeqno=\autoleqno}
\numeriadestra
\catcode`@=12
\def\as{\autosez}
\def\rs{\sref}
\def\ae{\autoeqno}
\def\re{\eqrefp}
\def\at{\autoteono}
\def\rt{\teoref}
\numerazionedoppia

\baselineskip 13pt

\def\rig#1{\smash{
   \mathop{\longrightarrow}\limits\sp {#1}}}

\def\dow#1{\Big\downarrow
   \rlap{$\vcenter{\hbox{$\scriptstyle#1$}}$}}

\font\bigg=cmbx12
\def\P#1{{\bf P}^#1}
\def\Z{{\bf Z}}
\def\pr{{\it Proof\ }}
\def\C{{\bf C}}

\def\Q{{\bf Q}}
\def\P{{\bf P}}
\def\O{{\cal O}}

\magnification=\magstep1

\riferimentifuturi
\def\S{{\cal S}}
\centerline{\bigg Unstable hyperplanes for Steiner bundles}
\centerline{\bigg and multidimensional matrices}

\centerline{\bf Vincenzo Ancona and Giorgio Ottaviani}
\footnote{\phantom{1}}{\rm
 Both authors were  supported by MURST and GNSAGA-INDAM}

\as{0} { Introduction}

  A multidimensional matrix of boundary format is an element $A\in
V_0\otimes\ldots\otimes V_p$
where
  $V_i$  is a  complex vector space of dimension $k_i+1$ for $i=0,\ldots ,p$
and
$$k_0=\sum_{i=1}^pk_i$$
We denote by $Det~A$ the hyperdeterminant of $A$ (see  [GKZ]).
 Let $e_0^{(j)},\ldots
,e_{k_j}^{(j)}$
be a basis in $V_j$ so that every $A\in V_0\otimes\ldots\otimes V_p$ has a
coordinate form
$$A=\sum a_{i_0,\ldots ,i_p}e_{i_0}^{(0)}\otimes\ldots\otimes e_{i_p}^{(p)}$$
 Let $x_0^{(j)},\ldots
,x_{k_j}^{(j)}$ be the coordinates  in $V_j$. Then $A$ has the
following different descriptions:
 \item {1)} A multilinear form
$$  \sum_{(i_0,\ldots ,i_p)} a_{i_0,\ldots
,i_p}x_{i_0}^{(0)}\otimes\ldots\otimes x_{i_p}^{(p)}$$
 \item {2)} An ordinary matrix $M_A = (m_{i_1i_0})$ of size $(k_1+1)\times
(k_0+1)$ whose entries are 
multilinear forms
$$m_{i_1 i_0} =  \sum_{(i_2,\ldots ,i_p)} a_{i_0,\ldots
,i_p}x_{i_2}^{(0)}\otimes\ldots\otimes x_{i_p}^{(p)}\ae{ma}$$
\item {3)} A sheaf  morphism  $f_A$ on the product $X=\P^{k_2}\times\ldots\times
\P^{k_p}$
$$\O_X^{k_0+1}\rig{f_A}\O_X(1,\ldots ,1)^{k_1+1} \ae{st}$$
The theorem 3.1
of chapter 14 of [GKZ] easily translates into:
\proclaim Theorem.
The following properties are equivalent
 \item {i)} $Det~A\neq 0$.
 \item {ii)} the matrix  $M_A$ has constant rank $k_1+1$ on $X=\P^{k_2}\times
 \ldots\times\P^{k_p}$.
\item {iii)} the  morphism  $f_A$  is  surjective so that $S^*_A = Ker f_A$
is a vector bundle of rank $k_0-k_1$.

The above remarks set up a basic link between non
degenerate multidimensional matrices of boundary
format and vector bundles on a product of projective spaces.
In the particular case $p=2$ the (dual) vector bundle
$S_A$  lives on the projective space
$\P^{n} , n=k_2$, and is  a  Steiner bundle as defined in  [DK]. We can
keep to $S_A$
 the name Steiner also for $p\ge3$.

The action of of $SL(V_0) \times\ldots\times SL(V_p)$ on $A$ translates to
an action on
the corresponding bundle in two steps:
first the action of $SL(V_0) \times  SL(V_1)$  leaves the bundle in the
same isomorphism class; then $SL(V_2) \times\ldots\times SL(V_p)$
 acts on the classes, i.e. on the moduli space of Steiner bundles.
It follows that the invariants of matrices for the action of
$SL(V_0) \times\ldots\times SL(V_p)$ coincide with the invariants of the
action of
$SL(V_2) \times\ldots\times SL(V_p)$
 on the moduli space of the corresponding bundles. Moreover the stable
points of both
actions correspond to each other.

The aim of this paper is to investigate the properties and the invariants
 of both the above actions.
When we look at  the vector bundles,  we restrict ourselves to the case $p=2$,
that is Steiner bundles on projective spaces. This is probably the first
case where Simpson's question
 ([Simp], pag. 11) about the natural $SL(n+1)$-action on the moduli spaces
of bundles  on
 $\P^n$ has been investigated.

The section 2 is devoted to the study of multidimensional matrices. We denote by
the same letter matrices in $V_0\otimes\ldots\otimes V_p$
and their projections in $\P (V_0\otimes\ldots\otimes V_p)$.
  In the theorem
\rt{t1}
we prove that a matrix $A\in\P (V_0\otimes\ldots\otimes V_p)$ of boundary format
such that
$Det~A\neq 0$ is not stable for the action of $SL(V_0) \times\ldots\times
SL(V_p)$
if and only if there is a coordinate system such that $a_{i_0\ldots i_p}=0$ for
$i_0>\sum_{t=1}^pi_t$. A matrix satisfying this condition is called
triangulable.
The other main results of this section are theorems
\rt{t2} and \rt{t3} which describe the behaviour of the  stabilizer subgroup
$Stab(A)$. In the remark \rt{discinv} we introduce a discrete $SL(V_0)\times 
SL(V_1)\times SL(V_2)$-invariant of nondegenerate matrices in $\P(V_0\otimes 
V_1\otimes V_2)$ and we show that it can assume only the values $0,\ldots ,
k_0+2,\infty$.
\medskip
The second part of the paper,  consisting  of   sections from $3$ to $6$,
  can be read independently of section 2, except that  we will use
theorem \rt{t1}
in two  crucial points (theorem \rt{3.2} and  section 6). In this part we 
study the Steiner bundles on
$\P^n=\P (V)$. As we  mentioned above, they are  rank-$n$
  vector bundles $S$ whose dual $S^*$ appears in an exact sequence
$$0\rig{} S^*\rig{} W\otimes\O\rig{f_A} I\otimes\O(1)\rig{} 0 \ae{1.0}$$
where $W$ and $I$ are complex vector spaces of
dimension $n+k$ and $k$ respectively.
The map $f_A$ corresponds to $A\in W^*\otimes
V\otimes I$ (which is of boundary format)
  and $f_A$ is surjective if and only if $Det~A\neq 0$.
 We denote   by $\S_{n,k}$ the family of Steiner
bundles described by a sequence as   \re{1.0}.  $\S_{n,1}$
contains only the quotient bundle. Important examples
of Steiner bundles are the Schwarzenberger bundles, whose  construction
goes back to
the pionieeristic
work of Schwarzenberger [Schw]. Other examples are the logarithmic
bundles $\Omega (log {\cal H})$ of meromorphic forms on $\P^n$
having at most logarithmic poles on a finite
union ${\cal H}$ of hyperplanes with normal crossing;
Dolgachev and Kapranov showed in [DK] that they are Steiner.
  The Schwarzenberger bundles  are  a special case of  logarithmic bundles,
when all the hyperplanes
osculate the same rational normal curve. Dolgachev and Kapranov proved a
Torelli type  theorem, namely that the logarithmic bundles are uniquely
determined up to
isomorphism by
the above union of hyperplanes, with a weak additional assumption. This
assumption was recently removed by Vall\`es[V2], who shares with us the idea of
looking at the scheme $W(S) = \{ H\in \P^{n\vee}|h^0(S^*_H)\neq 0\} \subset\P^{n\vee}$ of  unstable  hyperplanes of
a Steiner
bundle $S$.  
Vall\`es proves that  any $S\in\S_{n,k}$ with at least $n+k+2$ unstable
hyperplanes
with normal crossing  is a Schwarzenberger bundle and
 $W(S)$ is a rational normal curve.
 We strengthen this result by showing the
 following: for any $S\in\S_{n,k}$ any subset of closed points in
 $W(S)$ has always normal
 crossing (see the theor. \rt{1.3}). Moreover $S\in\S_{n,k}$ is logarithmic
if and only if
 $W(S)$ contains at least $n+k+1$ closed points (cor. \rt{Wlog1} and
\rt{Wlog2}).
In particular if $W(S)$
contains exactly $n+k+1$ closed points then
 $S\simeq\Omega(log~W(S))$. The Torelli
 theorem follows.

 It turns out
 that the length of $W(S)$ defines an interesting filtration into irreducibles
 subschemes
 of $\S_{n,k}$ which   gives also the discrete invariant of multidimensional
matrices of boundary format mentioned above. This filtration
 is well behaved   with respect to $PGL(n+1)$-action on $\P^n$  and also with
 respect to the classical notion of association reviewed in [DK]. Eisenbud and
Popescu realized in [EP] that the association is exactly what nowadays is
called Gale
transform. For Steiner bundles corresponding to $A\in W^*\otimes~V\otimes~I$
this operation amounts to exchange the role of $V$ with $I$, so that it
corresponds
to the transposition operator on multidimensional matrices.

 The Gale transform for Steiner bundles can be decribed by the natural
isomorphism

 $$\S_{n,k}/SL(n+1)\to \S_{k-1,n+1}/SL(k)$$

 Both quotients in the previous formula are isomorphic to the
 GIT-quotient $$\P(W^*\otimes V\otimes I)/SL(W)\times SL(V)\times SL(I)$$ which
 is a basic object in linear algebra.

 As an application of the tools developed in the first section we show that all
the points of $\S_{n,k}$ are semistable for the action
 of $SL(n+1)$ and we compute the stable points. Moreover we characterize 
 the Steiner bundles $S\in\S_{n,k}$ whose symmetry group (i.e. the group of
 linear projective transformations preserving  $S$) contains $SL(2)$ or
contains
 $\C^*$.

Finally  we mention  that $W(S)$ has a geometrical construction by
 means of the Segre variety. From this construction $W(S)$ can be easily
computed
 by means of current software systems. 

 We thank J. Vall\`es for   the useful discussions we had on the subject of this
paper.

\vfill\eject
\as{x}{Multidimensional Matrices of boundary Format and Geometric
Invariant Theory}

\proclaim Definition\at{x.1}. A $p+1$-dimensional matrix of boundary
format $A\in V_0\otimes\ldots\otimes V_p$
 is  called  triangulable if one of the following equivalent conditions holds:
\item{i)} there exist  bases  in $V_j$ such that $a_{i_0,\ldots ,i_p}=0$
for $i_0>\sum_{t=1}^p i_t$
\item{ii)} there exist  a vector space $U$ of dimension 2, a subgroup
$\C^*\subset SL(U)$ and isomorphisms $V_j\simeq S^{k_j}U$ such that  if
  $V_0\otimes\ldots\otimes V_p=
\oplus_{n\in\Z}W_n$ is the decomposition into direct sum of eigenspaces of
the induced representation, 
we have $A\in\oplus_{n\ge 0}W_n$

{\it Proof of the equivalence between i) and ii)}

 Let $x, y$ be a basis of $U$ such that $t\in\C^*$ acts on $x$ and $y$ as
 $tx$ and $t^{-1}y$. Set $e_k^{(j)}:=x^ky^{k_j-k}{{k_j}\choose k}\in
S^{k_j}U$ for $j>0$ and $e_k^{(0)}:=x^{k_0-k}y^k{{k_0}\choose k}\in S^{k_0}U$
so that $e_{i_0}^{(0)}\otimes\ldots\otimes e_{i_p}^{(p)}$ is a basis of
$S^{k_0}U\otimes\ldots\otimes S^{k_p}U$ which diagonalizes the action of $\C^*$.
  The  weight of $e_{i_0}^{(0)}\otimes\ldots\otimes e_{i_p}^{(p)}$
is $2\left( \sum_{t=1}^p i_t-i_0\right)$,  hence ii) implies i). The
converse is trivial.

The following definition agrees with the one in [WZ], pag. 639.

\proclaim Definition\at{x.2}. A $p+1$-dimensional matrix of boundary
format $A\in V_0\otimes\ldots\otimes V_p$   is called
diagonalizable if one of the following equivalent conditions holds
\item{i)}  there exist  bases in $V_j$ such that $a_{i_0,\ldots ,i_p}=0$
for $i_0\neq\sum_{t=1}^p i_t$
\item{ii)} there exist  a vector space $U$ of dimension 2, a subgroup
$\C^*\subset SL(U)$ and isomorphisms $V_j\simeq S^{k_j}U$ such that $A$ is
a fixed
point of
the induced action of $\C^*$.

The following definition agrees with the one in [WZ], pag. 639.

\proclaim Definition\at{x.3}.  A $p+1$-dimensional matrix of boundary format
 $A\in V_0\otimes\ldots\otimes V_p$ is  an identity if
one of the following equivalent
 conditions holds 
 \item{i)} there exist  bases in $V_j$ such that   
$$a_{i_0,\ldots ,i_p}=\left\{\matrix{0\quad\hbox{for}\quad
i_0\neq\sum_{t=1}^pi_t\cr
1\quad\hbox{for}\quad i_0=\sum_{t=1}^pi_t\cr}\right.$$
\item{ii)} there exist  a vector space $U$ of dimension $2$  and
isomorphisms $V_j\simeq S^{k_j}U$ such that
$A$ belongs to the unique one dimensional $SL(U)$-invariant subspace of
$S^{k_0}U\otimes S^{k_1}U\otimes\ldots\otimes S^{k_p}U$

The equivalence between i) and ii) follows easily from the following remark:
 the matrix  $A$ satisfies the condition ii) if and only if
it corresponds to the natural multiplication map
$S^{k_1}U\otimes\ldots\otimes S^{k_p}U\to S^{k_0}U$ (after a suitable
isomorphism
$U\simeq U^*$ has been fixed).

From now on, we consider the natural action  of
$SL(V_0)\times\ldots\times SL(V_p)$ on
$\P (V_0\otimes\ldots\otimes V_p)$.  We may suppose $p\ge 2$. 
The definitions of triangulable,
diagonalizable and identity
apply    to elements of  $\P (V_0\otimes\ldots\otimes V_p)$ as well. In
particular
all identity matrices fill a distinguished orbit in
$\P (V_0\otimes\ldots\otimes V_p)$. 
The hyperdeterminant of elements
of $V_0\otimes\ldots\otimes V_p$   was introduced by Gelfand,
Kapranov and Zelevinsky in [GKZ]. They  proved that the dual variety of the
Segre product
$\P(V_0)\times\ldots\times \P(V_p)$ is a hypersurface if and only if
$k_j\le\sum_{i\neq j} k_i$ for $j=0,\ldots ,p$ (which is  obviously true for a
matrix of boundary format). When the dual variety  is a hypersurface, its equation  is called the
hyperdeterminant of format
$(k_0+1)\times\ldots\times (k_p+1)$ and denoted by $Det$.
 The hyperdeterminant  
  is a homogeneous polynomial
function over $V_0\otimes\ldots\otimes V_p$ so that the condition $Det~A\neq 0$
is meaningful for $A\in \P (V_0\otimes\ldots\otimes V_p)$.
The function $Det$ is  $SL(V_0)\times\ldots\times SL(V_p)$-invariant, in
particular if $Det~A\neq 0$ then $A$ is semistable for the action of
$SL(V_0)\times\ldots\times SL(V_p)$.
We denote by $Stab~(A)\subset SL(V_0)\times\ldots\times SL(V_p)$
the  stabilizer subgroup of $A$ and by $Stab~(A)^0$ its connected component
containing
the identity. The main results of this section are the following.

\proclaim {Theorem \at{t1}}. Let $A\in \P (V_0\otimes\ldots\otimes V_p)$
 of boundary format such that $Det~A\neq 0$. Then
$$A\hbox{\ is triangulable}\iff A\hbox{\ is not stable for the action of\ }
SL(V_0)\times\ldots\times SL(V_p)$$

\proclaim {Theorem \at{t2}}. Let $A\in \P (V_0\otimes\ldots\otimes V_p)$
 of boundary format such that $Det~A\neq 0$.  Then
$$A\hbox{\ is diagonalizable}\iff \C^*\subset Stab(A)$$

We state the following theorem only in the case $p=2$, although we believe it
  true $\forall p\ge 2$. We  point out that in particular
$\dim~Stab~(A)\le 3$
which is a bound independent of $k_0, k_1, k_2$.

\proclaim{Theorem \at{t3}}. Let $A\in \P (V_0\otimes V_1\otimes V_2)$
of boundary format such that $Det~A\neq 0$.  Then there exists a $2$-dimensional
vector space $U$ such that $SL(U)$ acts over $V_i\simeq S^{k_i}U$ and
according to this action on $V_0\otimes V_1\otimes V_2$ we have
$Stab~(A)^0\subset SL(U)$. Moreover the following cases are possible
$$Stab~(A)^0\simeq\left\{\matrix{0&\cr \C\cr\C^*&\cr
SL(2)&\hbox{\ (this case occurs if and only if\ }A\hbox{\ is an
identity)}}\right.$$

\proclaim {Remark}. When $A$ is an  identity then $Stab~(A)\simeq SL(2)$.

Let $\Z_j$ be the finite set $\{0,\ldots ,j\}$. We set ${\cal B}:=
{\Z_{k_1}\times\ldots\times\Z_{k_p}}$.
A {\it slice (in the q-direction)} is  the subset  $\{ (\alpha_1,
\ldots ,\alpha_p)\in{\cal B} : \alpha_q=k\}$ for some
 $k\in\Z_q$. Two slices in the same direction are called {\it parallel}.
An {\it admissible path} is a finite sequence of elements  $(\alpha_1,
\ldots ,\alpha_p)\in{\cal B}$
starting from $(0,\ldots ,0)$, ending to $(k_1,\ldots ,k_p)$,
such that at each step exactly one $\alpha_i$ increases by $1$ and all other
remain equal. Note that each admissible path consists exactly of $k_0+1$
elements.

\proclaim Tom Thumb's lemma \at{ttl}. Put a mark (or a piece of bread)
on any element of any admissible path.
   Then  two
parallel slices contain the same number of marks.

\pr Any admissible path $P$  corresponds  to  a sequence of
$k_0$ integers  between
$1$ and $p$ such that the integer $i$ occurs exactly $k_i$ times. We call
this sequence the
code of the path $P$.  The occurrences
of the integer $i$ in the code divide all other integers different
from $i$ appearing in the code into $k_i+1$ strings (possibly empty); each
string
encodes the part of the path contained in one of the $k_i+1$ parallel slices.
The symmetric group
$\Sigma_{k_i+1}$ acts on the set ${\cal A}$ of all the admissible paths
by permuting the strings .  Let  $P^i_j$ the number of elements
(marks) of the path $P\in{\cal A}$ on the slice $\alpha_i=j$.
In particular $\forall\sigma\in\Sigma_{k_i+1}$ we have
$$\sum_{P\in{\cal A}}P^i_j=\sum_{P\in{\cal A}}\left(\sigma\cdot P\right)^i_j=
\sum_{P\in{\cal A}}P^i_{\sigma^{-1}(j)}$$
which proves our lemma.
\medskip

\proclaim \at{deg}. We  will often use the following well-known fact. If
$\ \O_X^k\rig{\phi} F$ is a
morphism of vector bundles on a variety $X$ with $k\le rank~F=f$ and
$c_{f-k+1}(F)\neq 0$  then the degeneracy locus $\{ x\in X|
rank(\phi_x)\le
k-1\}$   is nonempty of codimension $\le f-k+1$.

A square matrix with a zero left-lower submatrix with the NE-corner
on the diagonal has zero determinant. The following lemma generalizes this
remark
to multidimensional matrices of boundary format.

\proclaim {Lemma\at{lem1}}. Let $A\in V_0\otimes\ldots\otimes V_p$. Suppose
that in a suitable coordinate system there is $(\beta_1,\ldots ,\beta_p)\in{\cal B}$
such that
$a_{i_0\ldots i_p}=0$ for $i_k\le\beta_k$ ($k\ge 1$) and $i_0\ge\beta_0: =
\sum_{t=1}^p\beta_t$.
Then $Det~A=0$.

\pr The submatrix of $A$ given by elements  $a_{i_0\ldots i_p}$ satisfying
 $i_k\le\beta_k$ ($k\ge 1$)  
gives on $X=\P^{\beta_2}\times\ldots\times
\P^{\beta_p}$ the sheaf morphism
$$\O_X^{\beta_1+1}\to\O_X(1,\ldots ,1)^{\beta_0}$$ 
whose rank  by \rt{deg} drops  on a subvariety  of codimension
$\le  \beta_0-\beta_1=\sum_{t=2}^p\beta_t=\dim~
\P^{\beta_2}\times\ldots\times
\P^{\beta_p}$;
    hence there are nonzero vectors $v_i\in V_i^*$ for $1\le i\le p$ such that
$A(v_1\otimes\ldots\otimes v_p)=0$ and then $Det~A=0$ by the theorem 3.1
of chapter 14 of [GKZ].

\proclaim Lemma\at{lem2}. Let $p\ge 2$ and $a_j^i$ be integers
with $0\le i\le p$, $0\le j\le k_i$ satisfying
the inequalities
$a_j^0\ge a_{j+1}^0$ for $0\le j\le k_0-1\qquad$
$a_j^i\le a_{j+1}^i$ for $i>0\quad 0\le j\le k_i-1$
and the linear  equations
$$ \sum_{j=0}^{k_i}a_j^i=0\qquad\hbox{for\ }0\le i\le p$$
$$a_{\Sigma_{t=1}^p\beta_i}^0+a_{\beta_1}^1+\ldots +a_{\beta_p}^p=0
\quad\forall (\beta_1,\ldots ,\beta_p)\in{\cal B}$$
Then there is $N\in\Q$
such that
$$a_i^0=N(k_0-2i),\qquad a_i^j=N(-k_j+2i)\quad j>0$$
Moreover $N\in\Z$ if at least one $k_j$ is not even, and  $2N\in\Z$ if all
the $k_j$   are even.

\pr If $1\le s\le p$ and $\beta_s\ge 1$ we have the two equations
$$a_{\Sigma_{t=1}^p\beta_t}^0+a_{\beta_1}^1+\ldots +a_{\beta_s}^s+\ldots
 +a_{\beta_p}^p=0$$

$$a_{\Sigma_{t=1}^p\beta_t-1}^0+a_{\beta_1}^1+\ldots +a_{\beta_s-1}^s+\ldots
 +a_{\beta_p}^p=0$$
Subtracting  we obtain
$$a_{\Sigma_{t=1}^p\beta_t}^0-a_{\Sigma_{t=1}^p\beta_i-1}^0 =-
\left(a_{\beta_s}^s-a_{\beta_s-1}^s\right)$$
so that the right hand side does not depend on $s$.

Moreover  for $p\ge 2$  from the  equations

$$a_{\Sigma_{t=1}^p\beta_t}^0+a_{\beta_1}^1+\ldots +a_{\beta_q+1}^q\ldots +
a_{\beta_s-1}^s+\ldots +a_{\beta_p}^p=0$$

$$a_{\Sigma_{t=1}^p\beta_t}^0+a_{\beta_1}^1+\ldots +a_{\beta_q}^q\ldots +
a_{\beta_s}^s+\ldots +a_{\beta_p}^p=0$$

we get

$$ a_{\beta_q+1}^q-a_{\beta_q}^q=a_{\beta_s}^s-a_{\beta_s-1}^s $$
which implies   that the right hand side does not depend on
$s$ either. Let $a_{\beta_s}^s-a_{\beta_s-1}^s=2N\in\Z$.
Then  $a_t^s=a_0^s+2Nt$ for $t>0\quad s>0$.
By the assumption  $\sum_{t=0}^{k_s}a_t^s=0$
we get
$$(k_s+1)a_0^s+2N\sum_{t=1}^{k_s}t=0$$
that is
$$a_0^s=-k_sN$$
The formulas for $a_i^s$ and $a_i^0$ follow immediately.
If some $k_s$ is odd we have $2N\in\Z$ and $k_sN\in\Z$ so that
$N\in\Z$. $\diamond$

{\it Proof of the theorem \rt{t1}}.
 If $A$ is triangulable it is not stable.   
Conversely suppose $A$ not stable
and call again  $A$ a  representative  of $A$ in
$V_0\otimes\ldots\otimes V_p$.
By the Hilbert-Mumford criterion there exists a 1-dimensional parameter
subgroup
$\lambda\colon\C^*\to        
SL(V_0)\times\ldots\times SL(V_p)$ such that $\lim_{t\to 0}\lambda (t)A$ exists.
Let
$$a_0^s\le\ldots\le a_{k_s}^s\qquad 0\le s\le p $$
be the weights of the 1-dimensional parameter subgroup of  $SL(V_s)$
induced by  $\lambda$;
  with respect to a basis consisting of eigenvectors  the coordinate
$a_{i_0\ldots i_p}$
describes the eigenspace of  $\lambda$ whose  weight is
$a^0_{i_0}+a^1_{i_1}+\ldots +a^p_{i_p}$.
Recall that $$\sum_{j=0}^{k_i}a_s^i=0\qquad 0\le s\le p$$
We note that $\forall (\beta_1,\ldots ,\beta_k )\in{\cal B}$
we have
$$a_{\Sigma_{t=1}^p\beta_t}^0+a_{\beta_1}^1+\ldots +a_{\beta_p}^p\ge
0\ae{path} $$
otherwise the coefficient $a_{i_0\ldots i_p}$ is zero for
$i_k\le\beta_k\quad 1\le k\le p$ and $i_0\ge\sum_{t=1}^p\beta_t$
and  the lemma \rt{lem1} implies $Det~A=0$.
  The sum on  all $ (\beta_1,\ldots ,\beta_k )\in{\cal B}$ of  the left
hand side of \re{path}
 is nonnegative.
The contribution of $a^t$'s in this sum is zero by the Tom Thumb's
lemma \rt{ttl}.
Also the contribution of $a^0$'s is zero because it is zero on any
admissible path.
It follows that
$$a_{\Sigma_{t=1}^p\beta_t}^0+a_{\beta_1}^1+\ldots +a_{\beta_p}^p= 0\quad
\forall (\beta_1,\ldots ,\beta_k )\in{\cal B} $$
and by the lemma \rt{lem2} we get explicit expressions for the weights
which imply that $A$ is triangulable.
\medskip

{\it Proof of the theorem \rt{t2}}
We call again $A$ any representative  of $A$ in $V_0\otimes\ldots\otimes V_p$.
If $A$ is diagonal in a suitable basis $e_{i_0}^{(0)}\otimes\ldots\otimes
e_{i_p}^{(p)}$, we construct  a 1-dimensional parameter subgroup
$\lambda\colon\C^*\to
SL(V_0)\times\ldots\times SL(V_p)$
by the equation $\lambda (t) e_{i_0}^{(0)}\otimes\ldots\otimes
e_{i_p}^{(p)}:=t^{i_0-\sum_{t=1}^pi_t}
 e_{i_0}^{(0)}\otimes\ldots\otimes e_{i_p}^{(p)}$  so that $\C^*\subset
Stab(A)$.
Conversely let $\C^*\subset Stab(A)$. By the theorem \rt{t1} $A$ is
triangulable and by the lemma \rt{lem1} all diagonal elements
$a_{i_0\ldots i_p}$ with $i_0=\sum_{t=1}^pi_t$ are nonzero.
We can arrange the action on the representative  in order that the diagonal
corresponds to the zero eigenspace. Then the assumption
$\C^*\subset Stab(A)$  and the explicit expressions of the weights as in the
proof of the theorem \rt{t1} show that $A$ is diagonal. $\diamond$\medskip

We will prove theorem \rt{t3} by geometric arguments at the end of section 6.

\as{1} Preliminaries  about Steiner bundles

\proclaim Definition\at{defsteiner}. A Steiner bundle   over
$\P^n=\P (V)$
  is a vector bundle $S$ whose dual $S^*$ appears in an exact sequence
$$0\rig{} S^*\rig{} W\otimes\O\rig{f_A} I\otimes\O(1)\rig{} 0\ae{1.1}$$
where $W$ and $I$ are complex vector spaces of
dimension $n+k$ and $k$ respectively.

A Steiner bundle is  stable ([BS] theor. 2.7 or [AO] theor. 2.8) and  is
invariant by small
deformations
([DK] cor. 3.3). Hence the moduli space $\S_{n,k}$ of Steiner bundles
defined by    \re{1.1}
is isomorphic to an open subset of the  Maruyama moduli scheme of
stable bundles. On the other hand $\S_{n,k}$  is also isomorphic to the
GIT-quotient of a suitable
open subset of $\P (Hom(W,I\otimes V))$  for  the action of $SL(W)\times
SL(I)$
(see section \rs{4}). It is interesting to remark that these two
approaches give
two different compactifications of $\S_{n,k}$, but we do not pursue this
direction
in this paper. For other results about $\P (Hom(W,I\otimes V))$ see [EH].

\proclaim Definition \at{jump}.  Let  $S\in \S_{n,k}$ be a Steiner
bundle. A hyperplane
$H\in \P(V^*)$ is an  unstable hyperplane of $S$ if
$h^0(S^*_{|H})\neq 0$.    The set  $W(S)$ of the unstable hyperplanes  is
the  degeneracy
locus over  $\P(V^*)$ of the natural map  $H^1(S^*(-1))\otimes\O\to
H^1(S^*)\otimes\O (1)$, hence it
has a natural structure of scheme.
$W(S)$ is called the scheme of  the unstable
hyperplanes
of $S$. Note that since $h^0(S^*_{|H})\le 1$ ([V2]) the rank of the previous map 
drops at most by one.

\proclaim \at{wschema}. Let us describe more explicitly the map 
$H^1(S^*(-1))\otimes\O\to H^1(S^*)\otimes\O (1)$. 
From \re{1.1} it follows that $H^1(S^*(-1))\simeq I$ and 
$H^1(S^*)\simeq \left( V\otimes I\right) /W$. The projection 
 $V\otimes I\rig{B} \left( V\otimes I \right) /W$ can be interpreted as
a map $V\otimes H^1(S^*(-1))\rig{}H^1(S^*)$ which induces on $\P(V^*)$  the required 
morphism $H^1(S^*(-1))\otimes\O\to H^1(S^*)\otimes\O (1)$.

For a generic $S$,  $W(S)=\emptyset$.
Examples show that $W(S)$ can have a nonreduced structure.

We recall that if $D$ is a divisor with normal crossing then  $\Omega(log~D)$
is the bundle of meromorphic forms having at most logarithmic poles over $D$.
 If  ${\cal H}$ is the union of $m$ hyperplanes  $H_i$ with normal
crossing,
it is shown in [DK] that for  $m\le n+1$   $\Omega(log~{\cal H})$ splits
while for 
$m\ge n+2$ then $S\in\S_{n,k}$ where $k=m-n-1$.

The following is a simple consequence of [BS] theor. 2.5.

\proclaim Proposition \at{1.2}. Let $S\in \S_{n,k}$, then
$$h^0(S^* (t))=0\Longleftrightarrow t\le k-1$$

\pr $S^* (t)\simeq \wedge^{n-1}S(-k+t)$. The
$\wedge^{n-1}$-power of the sequence dual to \re{1.1} is
$$0\rig{}S^{n-1}I^*\otimes\O(-n+1-k+t)\rig{}S^{n-2}I^*\otimes W^*\otimes
\O(-n+2-k+t)\rig{}\ldots$$
$$\ldots\rig{}\wedge^{n-1}W^*\otimes\O(-k+t)\rig{}
\wedge^{n-1}S(-k+t)\rig{}0  $$
and from this sequence the result follows.

Let us  fix a basis in each of the vector spaces $W$ and $I$.
 Then  the morphism $f_A$ in \re{1.1} can be represented
by a $k\times (n+k)$ matrix $A$  (it was called $M_A$ in the introduction,
see \re{ma})
with entries  in $V$. In order to
simplify  the
notations we will use the same letter $A$ to denote also  its  class
in $\P (Hom(W,I\otimes V))$. $A$ has rank $k$
 at every point of $\P(V)$. Two such matrices represent  isomorphic
bundles if
and only if they lie in the same orbit of  the action of $GL(W)\times GL(I)$.

\proclaim \at{factsection}. In particular
$H^0(S^* (t))$
  identifies  to  the space  of $(n+k)\times 1$-column vectors $v$
with entries in $S^tV$ such that $$Av=0\ae{section}$$

Moreover $H\in W(S)$ (as closed point) if and only if there are nonzero
vectors $w_1$ of size $(n+k)\times 1$ and
$i_1$ of size $k\times 1$ both
with constant coefficients such that
$$Aw_1=i_1H\ae{1.2}$$

\proclaim \at{hyper}. According to the theorem stated in the introduction
 $A\in
Hom(W,V\otimes I)$ has nonzero hyperdeterminant if and only if it
corresponds to a
vector bundle. The locus in $\P(Hom(W,V\otimes I))$ where the hyperdeterminant
vanishes is an irreducible hypersurface of degree $k\cdot{{n+k}\choose k}$
([GKZ],chap. 14,
 cor.2.6).  It is interesting to remark that the above prop. \rt{1.2}
 can be proved also as a consequence of [GKZ], chap. 14, theor. 3.3.

 \proclaim \at{geomhyper}. The above description has a geometrical counterpart.
 Here $P(V)$ is the space of lines in $V$, dual
to the usual $\P$ of hyperplanes.
Consider in $ P(V\otimes I)$ the variety $X_r$ corresponding to elements
of $V\otimes I$ of rank $\le r$. In particular $X_1$ is the Segre variety
$P(V)\times P(I)$. Let $m=\min (n,k-1)$ so that $X_m$ is the variety
of non maximum rank elements.  Then  $A\in Hom(W,V\otimes I)$
 defines  a vector bundle if and only if it induces an embedding
$P(W)\subset P(V\otimes I)$  such that at every smooth point of  $X_m\cap P(W)$,
$P(W)$ and $X_m$ meet  transversally.
  This follows  from [GKZ], , chap. 14, prop. 3.14 and
chap. 1, prop. 4.11.

\proclaim \at{geom}. $W(S)$ has the following geometrical description.
Let $p_V$ be the  projection
of the Segre variety $P(V)\times P(I)$ on the $P(V)$.
Then 
$$W(S)_{red}=p_V[P(W)\cap \left( P(V)\times P(I)\right) ]_{red}$$
(according to the natural isomorphism
$P(V)=\P(V^*)$)
In fact  $i_1H$ in the formula \re{1.2} is a decomposable tensor in
$V\otimes I$.

\proclaim \at{geomschema}. About the scheme structure we remind that
$W(S)$ is the degeneration locus of the morphism
$I\otimes\O_{\P(V^*)}\rig{}{{V\otimes I}\over{W}}\otimes\O_{\P(V^*)}(1)$.
The following construction is standard.  The projective bundle
$\P=\P(I^*\otimes\O_{\P(V^*)})\rig{\pi}\P(V^*)$ is isomorphic to the Segre variety
$T=\P(V^*)\times\P(I^*)=P(V)\times P(I)$ and
$\O_{\P}(1)\simeq\O_T(0,1)$. The morphism
$$\C\rig{} {{V\otimes I}\over{W}}\otimes V^*\otimes I^*$$ 
defines
a section of $O_T(1,1)\otimes {{V\otimes I}\over{W}}$ with zero locus $Z=T\cap P(W)$.
Now assume that $dim~W(S)=0$, hence $dim~T=0$. 
By applying $\pi_*$ to the exact sequence
$$\O_T\otimes \left( {{V\otimes I}\over{W}}\right)^*\rig{}
\O_T(1,1)\rig{}\O_Z\rig{}0$$
 we get that the structure sheaf of $W(S)$
is contained in $\pi_*\O_Z$. We do not know if the 
equality always holds. In particular if $Z$ is reduced also $W(S)$ is reduced.
We will show in the Theorem \rt{doppio} that a multiple point occurs in $Z$ iff 
it occurs in $W(S)$.

\proclaim Theorem \at{1.3}. Let
 $S\in \S_{n,k}$ be a Steiner bundle. Then any set of distinct unstable hyperplanes 
of $S$ has normal crossing.

 \pr  We fix a coordinate system $x_0,\ldots ,x_n$ on $\P^n$ and a 
basis $e^1,\ldots e^{n+k}$ of  $W$.
Let $A$ be a matrix representing $S$.
If the thesis is not true, we may suppose that $W(S)$ contains the
hyperplanes  $x_0 = 0,\ldots
,x_j = 0,\sum_{i=0}^jx_i = 0$ for some $j$ such that $1\le j\le n-1$.  By
\re{1.2} there are $c^0\in W$, $b^0\in I$ such that $Ac^0=b^0x_0$. We
may suppose that the first coordinate of $c^0$ is nonzero, hence
$A\cdot\left[ c^0,e^2,\ldots ,e^{n+k}\right] =[b^0x_0,\ldots ]=A'$.

The matrix $A'$ still represents   $S$, hence by \re{1.2} there are
$c^1\in W$,
$b^1\in I$
such
that $A'c^1=b^1x_1$. At least one coordinate of $c^1$ after the first is
nonzero,
 say  the second. It follows that
$A'\cdot\left[e^1, c^1,e^3,\ldots ,e^{n+k}\right] = [b^0x_0, b^1x_1,\ldots
]=A''$ and again $A''$ represents $S$. Proceeding in  this way we get in
the end that
$$\left[b^0x_0,\ldots ,b^jx_j,\ldots\right]$$ is a matrix representing  $S$,
which   we denote
  again by $A$.

By \re{1.2}
there are $c=(c_1,\ldots ,c_{n+k})^t\in W$, $b\in I$ such that
$A\cdot c=\left[b^0x_0,\ldots b^jx_j,\ldots\right]\cdot c= b \sum_{i=0}^jx_i$

Now we distinguish two cases. If $c_i=0$ for $i\ge j+2$ we get
$b=c_1b^0=c_2b^1=\ldots =c_{j+1}b^j$, that is the submatrix of $A$ given
by the first
$j+1$ columns  has generically rank one. If we  take the $k\times
(n+k-j)$ matrix which has $b^j$ as first column and the last $n+k-j-1$ columns
of $A$ in the remaining places, we obtain a morphism
$$\O^k\rig{}\O\oplus\O(1)^{n+k-j-1}$$ which by \rt{deg} has rank $\le k-1$
on a nonempty
subscheme $Z$ of  $\P^n$.    It follows that also $A$ has rank $\le k-1$
on $Z$, contradicting
the assumption that $S$ is a bundle. So this case cannot occur.

In the second case  there exists a nonzero $c_i$ for some
$i\ge j+2$, we may suppose $c_{j+2}\neq 0$.
Then the matrix
$$A' = A\cdot\left[ e^1,\ldots ,e^{j+1},c,e^{j+3},\ldots ,e^{n+k}\right]=
\left[b^0x_0,\ldots b^jx_j,b\sum_{i=0}^jx_i\ldots\right]$$
  represents $S$.

 The last $n+k-j-2$ columns
of $A'$ define a sheaf  morphism $\O^k\rig{}\O(1)^{n+k-j-2}$ on the subspace
$\P^{n-j-1}=\{x_0=\ldots =x_j=0\}$
and again by \rt{deg}
 we find  a point
where the rank of $A$ is $\le k-1$. So  neither this case can occur.

 \proclaim Proposition\at{1.4}.
Let $S\in \S_{n,k}$ and let $\xi_1,\ldots ,\xi_s\in W(S)$, $s\le n+k$. 
There exists a matrix representing $S$  whose first $s$ columns are
$\left[ b^1\xi_1,\ldots , b^s\xi_s\right]$, where  the $b^i$ are vectors with
constant coefficients of size $k\times 1$. 
Moreover any $p$ columns
 among $b^1,\ldots ,b^s$ with $p\le k$ are independent. Conversely if the first
$s$ columns of a matrix representing $S$ have the form
$\left[ b^1\xi_1,\ldots , b^s\xi_s\right]$ then $\xi_1,\ldots ,\xi_s\in W(S)$.

 \pr The last assertion is obvious.  The proof of the existence of a
matrix $A$ representing $S$
having the required form is analogous to  that  of the theorem
\rt{1.3}. Then,  it is sufficient to prove that $b^1,\ldots ,b^p$ are
 independent.  Suppose
 $\sum_{i=1}^p b^i\lambda_i=0$. Let $\xi=\prod_{i=1}^p\xi_i$.
 Let  $c$  be the $(n+k)\times 1$ vector (whith coefficients  in
$S^{p-1}V$) whose
$i-th$ entry is
 $\lambda_i\xi/ \xi_i$ for $i=1,\ldots p$ and zero otherwise.
 It follows $A\cdot c=\xi\sum_{i=1}^p b^i\lambda_i=0$ and by $\re{section}$ we
  get a nonzero section   of $S^{*}(p-1)$, which contradicts the prop.
\rt{1.2}.

\proclaim \at{elm}. \it Elementary transformations.\rm
Consider $H=\{\xi=0\}\in W(S)$. The map $\O_H\to S^*_{|H}$ induces a
surjective map
$S\to\O_H$  and an exact sequence
$$0\rig{}S'\rig{}S\rig{}\O_H\rig{}0 \ae{1.2bis}$$
(see also [V2] theor. 2.1);
  it is easy to check (e.g. by Beilinson theorem) that
$S'\in\S_{n,k-1}$. According to [M] we say that $S'$ has been
obtained from $S$
by an elementary transformation. By the prop. \rt{1.4} there exists a matrix $A$ 
representing $S$ of the following form
$$A=\left[\matrix{\xi&*&\cdots&*\cr
0\cr
\vdots&&A'\cr
0\cr}\right]\ae{elmm}$$
where $A'$ is a matrix representing $S'$.
Since $h^0(S^*_{|H})\le 1$, $S'$ is uniquely determined by $S$ and $H$.

\proclaim Theorem \at{elmt}. With the above notations we have the 
inclusion of schemes
$W(S)\subset W(S')\cup H$. In particular
\item{i)} $length~W(S')\ge length~W(S)-1$
\item{ii)} If $dim~W(S')=0$ then $mult_HW(S')\ge  mult_HW(S)-1$, so that if 
$H$ is a multiple
point of $W(S)$, then $H \in W(S')$
\item{iii)} If $dim~W(S')=0$ then 
for any hyperplane $K\neq H\quad mult_KW(S')\ge mult_KW(S)$

\pr  The sequence  dual to \re{1.2bis}
$$0\rig{}S^*\rig{}S'^*\rig{}\O_H(1)\rig{}0 $$
gives the commutative diagram on $\P (V^*)$
$$\matrix{0&\rig{}&\O&\rig{}&H^1(S^*(-1))\otimes\O&\rig{}&H^1(S'^*(-1))\otimes\O&
\rig{}&0\cr
&&\dow{}&&\dow{}&&\dow{}\cr
0&\rig{}&H^0(\O_H(1))\otimes \O(1)&\rig{}&H^1(S^*)\otimes\O(1)&\rig{}&H^1(S'^*)
\otimes\O(1)&\rig{}&0}$$

It follows that the matrix $B'$ of the map
$$H^1(S'^*(-1))\otimes\O\to H^1(S'^*)\otimes\O (1)$$
 can be seen as a submatrix of the  matrix $B$ of the map
$$H^1(S^*(-1))\otimes\O\to H^1(S^*)\otimes\O (1)$$
  In a suitable system of coordinates:
$$B=\left[\matrix{
y_1&*\cr
\vdots&*\cr
y_n&*\cr
0&B'}\right]\ae{b}$$
where $(y_1,\ldots, y_n)$  is the ideal of $H$ (in the dual space).
It follows that
$$I(W(S'))\cdot (y_1,\ldots ,y_n)\subset I(W(S))$$
which concludes the proof.

\as{2} {The Schwarzenberger bundles}

 Let  $U$ be  a complex vector space of dimension $2$.
  The natural multiplication map $S^{k-1}U^*\otimes S^nU^* \to
S^{n+k-1}U^*$ induces the $SL(U)$-equivariant injective map
$S^{n+k-1}U\to S^{k-1}U\otimes S^nU$
and defines  a  Steiner bundle  on  $\P(S^nU)\simeq\P^n$  as  the dual of the
kernel of the surjective morphism
$$\O_{\P(S^nU)}\otimes S^{n+k-1}U\to \O_{\P(S^nU)}(1)\otimes S^{k-1}U$$
 It is called a Schwarzenberger bundle (see [ST],  [Schw]).
 Let us remark that in to the correspondence between Steiner bundles and 
multidimensional matrices mentioned in the introduction, the Schwarzenberger 
bundles correspond exactly to the identity matrices (see the def. \rt{x.3})

By interchanging the role of $ S^{k-1}U$ and  $S^nU$ we obtain also a
Schwarzenberger bundle
 on  $\P(S^{k-1}U)\simeq\P^{k-1}$  as  the
dual of the  kernel of the surjective morphism
$$\O_{\P(S^{k-1}U)}\otimes S^{n+k-1}U\to \O_{\P(S^{k-1}U)}(1)\otimes S^nU$$

Both  the above bundles are $SL(U)$-invariant.
We sketch the original Schwarzenberger construction for the first one.
The diagonal map $u\mapsto u^n$ and the isomorphism  $\P(S^nU)\simeq\P^n$
detect  a rational normal curve $\P(U)=C_n\subset\P^n$. In the same way  a
second
rational normal curve
$\P(U)=C_{n+k-1}$  arises in $\P(S^{n+k-1}U)$.   We  define  a morphism
$$\matrix{\P(S^nU)=S^n\P(U)&\to&Gr(\P^{n-1},\P(S^{n+k-1}U))\cr
n\ \hbox{points in\ }\P(U)&\mapsto &\hbox{Span of n points in\ }C_{n+k-1}}$$

The pullback of the dual of the universal bundle on the Grassmannian is a
Schwarzenberger bundle.

It is easy to check that if $S$ is a Schwarzenberger bundle then
$W(S)=C^*_n \subset \P(S^nU^*)$ (the dual rational normal curve).
(see e.g. [ST], [V1] ).

This can be explicitly seen    from the matrix form
given by [Schw, prop. 2]
$$M_A=\left[\matrix{ x_0&\ldots&x_n\cr
&\ddots&&\ddots\cr
&&x_0&\ldots&x_n }\right]\ae{schw}$$

Let $t_1,\ldots ,t_{n+k}$ be any distinct complex numbers.  Let $w$ be
  the $(n+k)\times (n+k)$
Vandermonde matrix whose
$(i,j)$ entry is $t_j^{(i-1)}$;  the $(i,j)$ entry of the product $M_Aw$ is
 $t_j^{(i-1)}\cdot (\sum_{k=0}^nx_kt_j^k)$; hence
  $\{ \sum_{k=0}^nx_kt^k=0\}\in W(S)$ $\forall t\in\C$ by the  prop. \rt{1.4}.
 On the other hand   $W(S)$ is $SL(U)$-invariant; if it where strictly
bigger than $C^*_n$ then it would
contain the hyperplane
$H =\{ x_0+x_1=0\}$,which lies in the next  $SL(U)$-orbit;  now   the equation
\re{1.2}
  implies   immediately $w_1 = i_1 =0$.

  In the Theorem \rt{nk2} we will need the following

\proclaim Lemma \at{SchwB}. Let $S$ be a Schwarzenberger bundle and let $(x_0,
\ldots ,x_n)$ be coordinates in $\P (V)$ such that $S$ is represented 
(with respect to suitable basis of $I$ and $W$) by the
matrix $M_A$ in \re{schw}. Let 
$(y_0, \ldots ,y_n)$ be dual coordinates in $\P (V^*)$. Then the morphism
$H^1(S^*(-1))\otimes\O\to H^1(S^*)\otimes\O (1)$
 (with respect to the obvious basis)  is represented by the matrix
$$B=\left[\matrix{y_1&-y_0&&&\cr
&y_1&-y_0&&\cr
&&\ddots&\ddots&\cr
&&&y_1&-y_0\cr
y_2&0&-y_0&&\cr
&\ddots&\ddots&\ddots&\cr
&&y_2&0&-y_0\cr
&&&y_2&-y_1\cr
y_3&0&0&-y_0&\cr
&\ddots&\ddots&\ddots&\ddots}\right]$$

\pr
By \rt{wschema} it is enough to check that the 
composition  $W\rig{A} V\otimes I\rig{B} \left( V\otimes I\right)/W$ is zero, which 
is straightforward.

\proclaim Theorem\at{2.moduli}. (Schwarzenberger, [Schw] theor. 1, see also
[DK] prop. 6.6).
The moduli space of Schwarzenberger bundles is $PGL(n+1)/SL(2)$, which is
the open subscheme
of the Hilbert scheme parametrizing rational normal curves.

In particular $W(S)$ uniquely determines  $S$ in the class of Schwarzenberger bundles.

\as{3} {A filtration  of $\S_{n,k}$ and the Gale transform of Steiner
bundles}

\proclaim Definition\at{3.5}.
$$\S_{n,k}^i:=\{S\in\S_{n,k}|length~W(S)\ge i\}$$

In particular $$\S_{n,k}=\S_{n,k}^0\supset\S_{n,k}^1\supset\ldots$$
We will see in a while (cor. \rt{3.00})
that $\S_{n,k}^{\infty}$ corresponds  to  Schwarzenberger bundles.

Each $\S_{n,k}^i$ is invariant for the action of   $SL(V)$   on
$\S_{n,k}$.
We will see in the section 6 that all the points of $\S_{n,k}$ are semistable
(in the sense of Mumford's GIT) for the action of $SL(V)$.

Let $\S$ be the open subset of $\P(Hom(W,V\otimes I))$ representing
 Steiner bundles.
The quotient $\S_{n,k}/SL(V)$ is isomorphic to $\S
/SL(W)\times SL(I)\times SL(V)$.

By interchanging  the role of $V$ and $I$,
also  $\S_{k-1,n+1}/SL(I)$ turns out to be isomorphic to $\S/SL(W)\times
SL(I)\times SL(V)$, so that
we obtain an isomorphism

$$\S_{n,k}/SL(n+1)\simeq\S_{k-1,n+1}/SL(k)$$

For any $E\in\S_{n,k}/SL(n+1)$ we will call the \it Gale transform \rm  of $E$
the corresponding class in $\S_{k-1,n+1}/SL(k)$ and we denote it by $E^G$.
 In [DK] the above construction is called association.
Here we follow  [EP].   Our Gale transform is a generalization of the one
in [EP]. In fact  in  the case
 $i=n+k+1$ Eisenbud
and Popescu in [EP] review the classical association between $PGL(n+1)$-classes
of $n+k+1$
 points of $\P^n$ in general position and $PGL(k)$-classes
of $n+k+1$
 points of $\P^{k-1}$ in general position and call it Gale transform. If we
take  the union
${\cal H}$ of $n+k+1$ hyperplanes with normal crossing in $\P^n$  (as
points in the dual projective space)
  the Gale transform (as points in the dual projective space) ${\cal H}^G$
consists of a $PGL(k)$-class of
$n+k+1$ hyperplanes with normal crossing in $\P^{k-1}$. As remarked in [DK]
  $[\Omega(log{\cal H})]^G\simeq [\Omega(log{\cal H}^G)]$. That is, the
Gale transform in our sense reduces to that
 in [EP] when the Steiner bundles are
logarithmic.  It is also clear that the $PGL$-class of
 Schwarzenberger bundles over $\P(V)$ corresponds  under the  Gale transform to
the $PGL$-class of Schwarzenberger bundles over $\P(I)$.

We point out  that one can define the Gale transform of a $PGL$-class of
 Steiner bundles but it is not possible to define the Gale transform of a
single
 Steiner
 bundle. This was implicit (but not properly written) in [DK]. Nevertheless
by a slight abuse we
will also speak about the Gale transform of a Steiner bundle
$S$, which will be any  Steiner bundle   in the class of the Gale transform
of $S~mod~SL(n+1)$.

 The following elegant
theorem due  to  Dolgachev and Kapranov  is a first beautiful
application of the Gale transform.

\proclaim Theorem\at{3.0}. (Dolgachev-Kapranov, [DK] theor. 6.8)
Any $S\in\S_{n,2}$ is a Schwarzen\-berger bundle.

\pr
$$\S_{n,2}/SL(n+1)\simeq\S_{1,n+1}/SL(2)$$
and it is obvious that a Steiner bundle on  the line $\P^1$ is
Schwarzenberger.

\proclaim Theorem\at{fund}. Two  Steiner bundles having in common
 ${n+k+1}$ distinct  unstable hyperplanes are isomorphic.

\pr We prove that if $S$  is a Steiner bundle  such that the   hyperplanes
$\{\xi_i=0\}$  for
$i=1,\ldots ,n+k+1$ belong to $W(S)$, then $S$ is uniquely determined. By
the prop. \rt{1.4} there exist
column vectors
 $a_i\in\C^k$    such that $S$ is represented by the matrix
$\left[ a^1\xi_1,\ldots ,a^{n+k}\xi_{n+k}\right]$
Moreover by \re{1.2} there are $b\in\C^{n+k}$ and
$c\in\C^k$ such that
 $$\left[ a^1\xi_1,\ldots , a^{n+k}\xi_{n+k}\right] b=c\xi_{n+k+1}$$
We claim that {\it all} the components of $b$ are nonzero.  The last
formula can be written

$$\left[ a^1b^1,\ldots,a^{n+k}b^{n+k},-c\right]\cdot\left[ \xi_1,\ldots
,\xi_{n+k+1}\right] ^t =0$$

where in the right matrix we identify $\xi_i$  with  the $(n+1)\times 1$ vector
given by the coordinates of the corresponding hyperplane.
We may suppose that there exists $s$ with $1\le s\le n+k-1$ such that
$b_i=0$ for $1\le i\le s$ and $b_i\neq 0$ for $s+1\le i\le n+k$.
If $s\ge k$, it follows that  $n+1$ hyperplanes among the $\xi_i$ have a
nonzero syzygy,
which contradicts  the prop. \rt{1.3}. Hence $s\le k-1$ and we have
$$\left[ a^{s+1}b^{s+1},\ldots,a^{n+k}b^{n+k},-c\right]\cdot\left[
\xi_{s+1},\ldots
,\xi_{n+k+1}\right] ^t =0$$
The rank of the right matrix is $n+1$, hence the rank of the left matrix
is $\le k-s$, in particular the first $k-s+1$ columns are dependent and
this contradicts
the prop. \rt{1.4}. This proves the claim.

In particular
$\left[ a^1,\ldots ,a^{n+k},-c\right]\cdot B=0$ where
$$B=Diag(b_1, \ldots ,b_{n+k},1)\cdot \left[ \xi_{1},\ldots
,\xi_{n+k+1}\right] ^t$$
 is a $(n+k+1)\times
(n+1)$ matrix with constant entries of  rank $(n+1)$. Therefore  the matrix
$\left[ a^1,\ldots a^{n+k},-c\right]$
is uniquely determined up to the  (left) $GL(k)$-action, which implies
that $S$ is uniquely determined up to isomorphism.

\proclaim Corollary\at{logarithmic}.  A Steiner bundle is logarithmic if
and only if it  admits at least
$(n+k+1)$ unstable hyperplanes.
 
\pr In fact ${\cal H}\subset W(\Omega(log{\cal H}))$ by the formula (3.5) of [DK] and
the prop. \rt{1.4}. 

\proclaim Corollary\at{3.00}. (Vall\`es, [V2] theor. 3.1)
 A Steiner bundle is Schwarzenberger if and only if it admits  at least
$(n+k+2)$ unstable hyperplanes.
In particular $\S_{n,k}^{\infty}$ coincides with the moduli space of
Schwarzenberger bundles.

\pr Let $S$ be a Steiner bundle, and   $H\in W(S)$.  Let us consider
the elementary transformation \rt{elm}
$$0\rig{}S'\rig{}S\rig{}\O_H\rig{}0$$
where
$S'\in\S_{n,k-1}$;  by the theorem \rt{elmt} has $n+k+1$ unstable  hyperplanes.
  Picking  $H'\in W(S')$ and repeating the above procedure
  after $(k-2)$ steps  we reach a  $S^{(k-2)}\in\S_{n,2}$;
by the theor. \rt{3.0} $S^{(k-2)}$ is a Schwarzenberger bundle. In
particular the
remaining $n+4$ unstable  hyperplanes lie on a rational normal curve.   It
is then clear that
  any   subset of $n+4$ hyperplanes in $W(S)$ lies on a
rational normal curve. Since there is a unique rational normal curve
through $n+3$
points in general position, it follows  that $W(S)$ is contained in a rational
normal curve, so that  $S$ is a Schwarzenberger bundle by the theorem \rt{fund}.

\proclaim Theorem\at{3.1}. Let $n\ge 2$, $k\ge 3$.
\item{i)} $\S_{n,k}^i$ for $0\le i\le n+k+1$ is an irreducible unirational
closed subvariety of  $\S_{n,k}$ of dimension
$(k-1)(n-1)(k+n+1)-i[(n-1)(k-2)-1]$
\item{ii)}$\S_{n,k}^{n+k+1}$ contains as an open dense subset
 the variety of Steiner logarithmic bundles
which coincides with the open subvariety of $Sym^{n+k+1}\P^{n\vee}$
consisting of
hyperplanes in $\P^n$ with normal crossing.

\pr
 (ii) follows from the  theorem  \rt{fund}. The irreducibility
 in (i) follows from the geometric construction \rt{geom}. The numerical
computation
 in (i) is performed (for $i\le n+k$) by adding $i(n+k-1)$ (moduli of $i$
 points in $\P(V)\otimes\P(I)$)
to $n(k-1)(n+k-i)$ (dimension of Grassmannian of linear $\P^{n+k-1}$ in
$\P(V\otimes I)$
 containing the span of the above $i$ points) and subtracting $k^2-1$
($\dim~SL(I)$).

\proclaim Remark \at{nk23}. In the case $(n,k)=(2,3)$ the generic Steiner
bundle is logarithmic
(this  was  remarked in [DK], 3.18).
In fact the generic $\P^4$ linearly embedded in $\P^8$ meets the Segre variety
$\P^2\times\P^2$ in $deg~\P^2\times\P^2=6=n+k+1$ points.

\proclaim Remark. The dimension of  $\S_{n,k}^i/SL(n+1)$   is equal to
$(n+k+1-i)[(k-2)(n-1)-1]+n(k-1)$ for $k\ge 3$, $n\ge 2$, $0\le i\le n+k+1$
and it is $0$ for $i\ge n+k+2$.

\proclaim \at{consegre}. The cor. \rt{3.00} implies the following property
of    the Segre
variety:  if a generic linear $P(W)$ meets $P(V)\times P(I)$ in $n+k+2$
points then meets it in infinitely many points.

 \proclaim Theorem\at{3.2}. Consider a nontrivial (linear) action
of $SL(2)=SL(U)$ over $\P^n$. If a Steiner bundle  is $SL(2)$-invariant
 then it is a Schwarzenberger bundle and $SL(U)$ acts over $\P^n=\P(S^nU)$.
Hence $\S_{n,k}^\infty$
 is the   subset of the fixed points of the action of $SL(2)$
 on $\S_{n,k}$.

 \pr By the theorem \rt{t1} there exists a coordinate system such that all
the entries (except the first)
 of  the first column of  the matrix
representing the Steiner bundle $S$ are zero. By the prop. \rt{1.4}
$W(S)$ is nonempty.
 By the assumption $W(S)$ is $SL(2)$-invariant and closed, it
follows that $W(S)$ is a union of rational curves and of simple points.
If $W(S)$ is infinite we can apply the corollary \rt{3.00}. If
$W(S)$ is finite we argue by induction on $k$ .We pick up $H\in W(S)$ 
and we consider the elementary transformation
$0\rig{}S'\rig{}S\rig{}\O_H\rig{}0$.
We get $\forall g\in SL(U)$ the diagram
$$\matrix{S&\rig{\phi}&\O_H\cr
\dow{i}\cr
g^*S&\rig{g^*\phi}&\O_H}$$
Since $h^0(S^*_{|H})\le 1$ we get that $\phi$ and $g^*\phi\circ i$ coincide
up to a scalar multiple. We obtain a commutative diagram
$$\matrix{0&\rig{}&S'&\rig{}&S&\rig{}&\O_H&\rig{}&0\cr
&&\dow{}&&\dow{\simeq}&&\dow{\simeq}\cr
0&\rig{}&g^*S'&\rig{}&g^*S&\rig{}&\O_H&\rig{}&0\cr}$$
It follows that $S'\simeq g^*S'$, hence $SL(U)\subset Sym(S')$ and by the 
inductive assumption $S'$ is Schwarzenberger and $SL(U)$ acts over $\P^n=\P(S^nU)$.
Hence $W(S)$ is infinite and we apply again the corollary \rt{3.00}.

\proclaim Corollary\at{Wlog2}.    If ${\cal H}$ is the
 union of $n+k+1$
 hyperplanes with normal crossing then

$$W(\Omega(log~{\cal H}))=\left\{\matrix{ {\cal H} &\hbox{
 when\ }{\cal H}\hbox{\   does not osculate a rat. normal curve}\cr
C_n&\hbox{
 when\ }{\cal H}\hbox{\ osculate the rat. normal curve\ } C_n\cr
&\hbox{(this case occurs iff\ } \Omega(log~{\cal H}) \hbox{\ is Schwarzenberger)}}
\right.$$

\pr  $\cal H \subset $$\Omega(log~{\cal H}) )$ by the prop. \rt{1.4}.
The result
follows by the theorem \rt{fund} and the cor. \rt{3.00}.

\proclaim Corollary\at{Wlog1}.  Let  $S\in\S_{n,k}$ be a Steiner bundle.
 If  $W(S)$ contains at least $n+k+1$ hyperplanes then
  for every subset ${\cal H}\subset W(S)$ consisting of
$n+k+1$ hyperplanes  $S\simeq\Omega(log~{\cal H})$, in particular  $S$ is logarithmic.

\proclaim Corollary \at{torelli}. (Torelli theorem, see [DK] for $k\ge n+2$
or [V2] in
general).
Let ${\cal H}$ and ${\cal H}'$ be two finite unions of $n+k+1$ hyperplanes
with normal crossing
in $\P (V)$ with $k\ge 3$ not osculating any rational normal curve. Then
$${\cal H}={\cal H}'\Longleftrightarrow \Omega(log~{\cal H})\simeq
\Omega(log~{\cal H}')$$

\proclaim Theorem \at{nk2}. Let $S\in\S_{n,k}$ be a Steiner bundle. If 
$length~W(S)\ge n+k+2$
then $length~W(S)=\infty$ and $S$ is Schwarzenberger.

\pr We proceed by induction on $k$. If $k=2$ the result follows from
the theorem \rt{3.0}, so we can suppose $k \ge 3$.
Let us pick any  $H\in W(S)$ and perform the elementary
transformation \re{1.2bis}. Then $S'\in\S_{n,k-1}$ and by the theor. \rt{elmt}
i)  $length~W(S')\ge n+k+1$, so that
by induction $S'$ is Schwarzenberger, in particular $W(S')$ is a rational
normal curve $C_n$.

It follows that $S$ is represented by the matrix
$$M_A=\left[\matrix{x_0&f_1&f_2&\ldots&&&f_{n+k-1}\cr
&x_0&x_1&\ldots&x_n\cr
&&\ddots&\ddots&&\ddots\cr
&&&x_0&x_1&\ldots&x_n}\right]$$
where
$f_i=-\sum_{j=1}^nc^i_jx_j$.
It is easy to check after the lemma \rt{SchwB} (and the proof of Theor. \rt{elmt})
that the morphism $H^1(S^*(-1))\otimes\O\to H^1(S^*)\otimes\O (1)$
is represented by the matrix
$$B=\left[\matrix{
y_1&c^1_1y_0&c^2_1y_0&\ldots&c^{k-2}_1y_0&\sum_{h=0}^nc^{k+h-1}_1y_h\cr
y_2&c^1_2y_0&c^2_2y_0&\ldots&c^{k-2}_2y_0&\sum_{h=0}^nc^{k+h-1}_2y_h\cr
\vdots&\vdots&\vdots&&\vdots&\vdots\cr
y_n&c^1_ny_0&c^2_ny_0&\ldots&c^{k-2}_ny_0&\sum_{h=0}^nc^{k+h-1}_ny_h\cr
&y_1&-y_0\cr
&&\ddots&\ddots\cr
\cr
&&&&y_1&-y_0\cr
&y_2&0&-y_0\cr
&&\ddots&\ddots&\ddots\cr
&&&y_2&0&-y_0\cr
&&&&y_2&-y_1\cr
&y_3&0&0&-y_0\cr
&&\ddots&\ddots&\ddots&\ddots}\right]$$
By the Theorem \rt{elmt} we have that 
$$length~\left( W(S)\cap C_n\right) \ge n+k+1\ae{wcn}$$
The points of $C_n$ are parametrized by
$y_i=t^i$ and $W(S)\cap C_n$ is given by the $k\times k$ minors of $B$
where we substitute $y_i=t^i$.
It is sufficient to look at the first $n+k-2$ rows because the others are linear 
combination of these. The first two rows and the last $k-2$ give the submatrix
$$\left[\matrix{
t&c^1_1&c^2_1&\ldots&c^{k-2}_1&\sum_{h=0}^nc^{k+h-1}_1t^h\cr
t^2&c^1_2&c^2_2&\ldots&c^{k-2}_2&\sum_{h=0}^nc^{k+h-1}_2t^h\cr
&t&-1\cr
&&t&-1\cr
&&&\ddots&\ddots\cr
&&&&t&-1}\right]$$
whose determinant is given up to sign by
$$t^{n+k}c^{n+k-1}_1+t^{n+k-1}(c^{n+k-2}_1-c^{n+k-1}_2)+\ldots
+t^2(c^1_1-c^2_2)-tc^1_2;\ae{pol}$$
by \re{wcn} all the coefficients of this polynomial are zero.
When $n=2$ this is enough to conclude that $M_A$ represents a Schwarzenberger bundle
because the matrix $M_A$ reduces to \re{schw} after a Gaussian elimination on the 
rows.
If $n\ge 3$ we have to look also at other minors.
For example the minor given by the first, third and the last $k-2$ rows is
$$\left[\matrix{
t&c^1_1&c^2_1&\ldots&c^{k-2}_1&\sum_{h=0}^nc^{k+h-1}_1t^h\cr
t^3&c^1_3&c^2_3&\ldots&c^{k-2}_3&\sum_{h=0}^nc^{k+h-1}_3t^h\cr
&t&-1\cr
&&t&-1\cr
&&&\ddots&\ddots\cr
&&&&t&-1}\right]$$
whose determinant is equal up to sign to
$$t^{n+k+1}c^{n+k-1}_1+t^{n+k}c^{n+k-2}_1+t^{n+k-1}(c^{n+k-3}_1-c^{n+k-1}_3)+\ldots
+t^3(c^1_1-c^3_3)-t^2c^2_3-tc^1_3$$
By \re{pol} the leading term $c^{n+k-1}_1$ vanishes and the degree drops
so that by \re{wcn} also the coefficients of this last polynomial vanish.
The reader can convince himself that the same argument of the case $n=2$ works 
also in this case.

  We remark that the above proof does not use the
cor. \rt{3.00} and gives a second proof of this corollary.

\proclaim Remark. There are examples of Steiner bundles $S\in\S_{n,k}$ such that
$length~W(S)=n+k+1$  and $W(S)$, as a set,  consists of only one point.

\proclaim Remark \at{discinv}. The above theorem shows that the only possible 
values for $length~W(S)$ are $0,1,\ldots ,n+k+1,\infty$. With the notations of 
section 2, every 
multidimensional matrix $A\in V_0\otimes V_1\otimes V_2$ of boundary format such 
that $Det~A\neq 0$ has a $GL(V_0)\times GL(V_1)\times GL(V_2)$-invariant
$$w(A):=length~W(\ker~f_A)^*$$
which can assume only the values $0,1,\ldots , dim V_0+1,\infty$.

\as{4} {Moduli spaces  of Steiner bundles and Geometric Invariant Theory}

Let ${\cal S}\subset \P(Hom(W,V\otimes I))$ be the open subset consisting
of $\phi\colon W\to V\times I$ such that for every nonzero $v^*\in
V^*$ the composite 
  $v^*\circ\phi\colon W\to I$ has maximum rank. By \rt{hyper},
${\cal S}$ is the complement of a hypersurface, and it is
  invariant for the natural action of $SL(W)\times SL(I)$. By interchanging
the roles of
$V$ and
$I$
  (or, in the language of the previous section, by performing the Gale 
transform) 
it is easy to check that  ${\cal S}$ coincides with the open subset of
$\phi\colon W\to V\times I$ such that for every nonzero $i^*\in
I^*$ the composite 
  $i^*\circ\phi\colon W\to V$ has maximum rank.

\proclaim Lemma \at{4.1}.
Every point of ${\cal S}$  is stable for the action of
$SL(W)\times SL(I)$.

\pr. Suppose that $A\in {\cal S}$ is not stable. Then by the Hilbert-Mumford
criterion
there exists a one-parameter subgroup
$\lambda (t)\colon\C^*\to
SL(W)\times SL(I)$ such that $\lim_{t\to 0} \lambda (t) A$ exists.
We may suppose that the two projections of $\lambda(t)$ on the factors
act diagonally with weights $\beta_1\le\beta_2\le\ldots \le\beta_k$
and $\gamma_1\le\gamma_2\le\ldots\le\gamma_{n+k}$ such that
$\sum_i\beta_i=\sum_j\gamma_j=0$.

We claim that there exists $p$ such that $1\le p\le k$ and
$\beta_p+\gamma_{k+1-p}<0$. Otherwise we get
$0\le\sum_{i=1}^k(\beta_i+\gamma_{k+1-i})=\sum_{i=1}^k\beta_i+\sum_{j=1}^k\gamma_j=
\sum_{j=1}^k\gamma_j\le k\gamma_k$, hence $0\le\gamma_k$.
If $\gamma_{n+k}>0$ we have $0\le\sum_{j=1}^k\gamma_j<
\sum_{j=1}^{n+k}\gamma_j=0$
which is a contradiction. If $\gamma_{n+k}=0$ then $\gamma_j=0$ $\forall j$ and
the claim is obvious. It follows that $\beta_i+\gamma_j< 0$ for
$i\le p$ and $j\le k+1-p$. Hence the first $p\times (k+1-p)$ block of the matrix
corresponding to $A$ is zero. The first $p$ rows of $A$ define a morphism
$\O^p\to\O(1)^{n+p-1}$ that by \rt{deg} drops rank in codim $\le n$,
contradicting the
fact
that  $A$  has maximum rank at every point.$\diamond $

\proclaim Theorem \at{4.3}.
Every point of ${\cal S}$  is semistable for the action of
$SL(W)\times SL(V)\times SL(I)$.

\pr. By \rt{hyper},  ${\cal S}$ is the complement of a
$SL(W)\times SL(V)\times SL(I)$-invariant
hypersurface ([GKZ], chap. 14, prop. 1.4).

\proclaim Corollary\at{4.4}. Every point of $\S_{n,k}$ is semistable for the action
of $SL(V)$ (with respect  to the natural polarization of $\S_{n,k}$ as
GIT-quotient).

\pr We look at the hyperdeterminant as a polynomial in the coordinate
ring of the
GIT quotient
$\P(Hom(W,V\otimes I))/SL(W)\times SL(I)\supset\S_{n,k}$ which is invariant
by the action of $SL(V)$.

\proclaim Theorem\at{rt1}.  An element  $A\in\S_{n,k}$
is not stable for  the action of $SL(n+1) = SL(V)$ if and only if there is a 
coordinate system such 
that the ordinary matrix $M_A$ (with entries in $V$) associated to  $A$
(see \re{ma}) has the
triangular form
$M_A=\sum_{j=0}^nA^mx_m$,  where the $(i,j)$-entry $a_{ij}^m$ of $A^m$ is zero for
$j<i+m$

\pr It is a reformulation of  the theorem \rt{t1} in the case $p=2$.

\proclaim Proposition\at{doppio}. Let $S$ be a Steiner bundle. The following properties are equivalent:
\item {i)} there is a   hyperplane $H$  which  is a multiple point for 
$W(S)$,  or $S$ is a Schwarzenberger bundle;
\item {ii)} there is a coordinate system such that $H = \{x_0=0\}$ and the matrix 
$M_A = \sum_{j=0}^nA^mx_m$  satisfies $a_{ij}^0=0$ for $j<i$, $j=1,2$    and
$a_{ij}^m=0$ for $m\ge 1$, $j\le i $, $j=1,2$ .

\pr  By \rt{geomhyper}, \rt{geom} and \rt{geomschema} (with the same  notations)
  if the condition i) occurs then $S$ is Schwarzenberger or
$Z$ has a multiple point. In both cases  there is  some point of $P(V)\times P(I)$
whose tangent space intersects $P(W)$  
 in a subspace of positive dimension.  The tangent 
space at a point $[v_0\otimes i_0]\in P(V)\times P(I)$ is the 
span of the two linear subspaces $P(V\otimes <i_0>)$ and $P(<v_0>\otimes I)$,
so that any point of the tangent space has the form $[v_1\otimes i_0+v_0\otimes i_1]$.
If the point $[v_1\otimes i_0+v_0\otimes i_1]$ with $v_0\neq v_1$, $i_0\neq i_1$
belongs to $P(W)$ it is easy to 
check that the matrix of $S$ satisfies ii). Conversely if the matrix of $S$
satisfies ii) then according to \re{elmm}
we can perform twice the elementary transformation at the hyperplane $H$
corresponding to $v_0$. Let $y_0,\ldots ,y_n$ be coordinates in $\P(V^*)$ such that
the ideal of $\{ H\}$ is defined by $y_1,\ldots ,y_n$.
The matrix $B$ in \re{b} has the form
$$B=\left[\matrix{y_1&g_1(y_0,\ldots ,y_n)&*\cr
\vdots&\vdots&*\cr
y_n&g_n(y_0,\ldots ,y_n)&*\cr
0&y_1&*\cr
\vdots&\vdots&*\cr
0&y_n&*\cr
0&0&B'}\right]$$
where $g_i$ are linear forms.
It is straightforward to check that the maximal minors of the restriction of $B$
to the line parametrized by
$y_0=1$, $y_i=tg_i(1,0,\ldots ,0)$ for $i=1,\ldots ,n$ 
have a multiple root for $t=0$, hence either $H$ is a multiple point of $W(S)$
or $W(S)$ is a curve and $S$ is Schwarzenberger by the Cor. \rt{3.00}.
$\diamond$

\proclaim Corollary \at{doppiobis}. With the notations of \re{1.2bis}
if $H$ is a multiple 
point of $W(S)$ then $H\in W(S')$.

\pr By the theorem \rt{doppio} the matrix $A$ representing $S$ has the form
\re{elmm} where $A'$ has the same form.

\proclaim Corollary \at{4.6}. If $S\in\S_{n,k}$ is not stable for the action
of $SL(V)$
then $S\in\S_{n,k}^2$

\pr From theor. \rt{rt1} and prop. \rt{doppio}.

\proclaim Remark. We conjecture  that if $S\in\S_{n,k}$ ($k\ge 3, (n,k)\neq (2,3)$)
is not stable for the action
of $SL(V)$ then $S\in S_{n,k}^3$  
and moreover  $S$ is Schwarzenberger or $W(S)$ has a point of multiplicity 
at least $3$.
We can prove that $S$ is Schwarzenberger or,
 in the notations of \rt{geom}, $Z=P(W)\cap \left( P(V)\times P(I)\right)$ has a 
point of multiplicity at least $3$.

 \proclaim Theorem\at{4.7}.  Let $S\in\S_{n,k}$ be a Steiner bundle. The following two
conditions are equivalent
\item{i)} $Sym(S)\supset\C^*$
\item{ii)} there is a coordinate system such that the matrix of $S$ has the 
diagonal form
$$\left[\matrix{ a_{0,1}x_0&\ldots&a_{n,1}x_n\cr
&\ddots&&\ddots\cr
&&a_{0,k}x_0&\ldots&a_{n,k}x_n }\right]$$

\pr It is a reformulation of the theorem \rt{t2} in the case $p=2$.  

\proclaim Corollary\at{cstar}. Let $S\in\S_{n,k}$ be a Steiner bundle  such that $Sym(S)\supset\C^*$.
 Then the $\C^*$-action on $\P^n$ has exactly $n+1$ fixed points whose weights are 
 proportional to $-n,-n+2,\ldots ,n-2, n$.

\pr The statement follows from the def. \rt{x.2}.

\proclaim Corollary\at{cstar1}.  Let $S\in\S_{n,k}$  be a Steiner bundle such that $Sym(S)\supset\C^*$.
Then either $W(S)$ is a rational normal curve and $S$ is a Schwarzenberger bundle,  or
$W(S)$ has only two closed points, namely   the two fixed points 
of the dual $\C^*$-action on $\P^{n\vee}$  having minimum and maximum weights.

\pr If $S$ is not  Schwarzenberger, $W(S)$ is finite   (by the cor. \rt{3.00}); since it is 
$Sym(S)$-invariant,  
  it must be contained 
in the $n+1$ fixed points of the $\C^*$-action on $\P^{n\vee}$. It is now easy to check, with the 
notations of \re{1.2}, that
the equation
$$\left[\matrix{ a_{0,1}x_0&\ldots&a_{n,1}x_n\cr
&\ddots&&\ddots\cr
&&a_{0,k}x_0&\ldots&a_{n,k}x_n }\right]\cdot w_1=i_1\cdot x_j$$
has nonzero solutions only for $j=0,n$.

\proclaim Proposition\at{4.8}. A logarithmic bundle  in $\S_{n,k}$
 which is not    stable for the action
of $SL(n+1)$  is  Schwarzenberger.

\pr The proof is by induction on $k$. For $k=2$ the result is true by the theorem 
\rt{3.0}. By the theorem \rt{rt1} there exists a triangular matrix corresponding 
to $S$.  Then  $H=\{ x_0 = 0\}$ is an unstable hyperplane of $S$. By  \rt{elm} 
there is an elementary transformation
$$0\rig{}S'\rig{}S\rig{}\O_H\rig{}0$$
where also  $S'$ is logarithmic (by the theor. \rt{elmt} and the coroll. \rt{Wlog2}).
Moreover the matrix representing $S'$ is again triangular by $\re{elmm}$.
$S'$ is not stable by the theor. \rt{rt1}. By induction  
$S'$ is Schwarzenberger and $  W(S')=C_n$ is a rational normal curve. For every $K\in W(S)$, $K\neq H$, we have 
$K\in W(S')=C_n$ by the theor. \rt{elmt}.
  The crucial point is that in this case also $H\in W(S')=C_n$; this can be checked by 
looking at the matrix of $S'$.
  Hence every closed point of  $W(S)$ 
lies in $C_n$ and by the  theorem \rt{fund} $S$ is isomorphic to the 
Schwarzenberger bundle determined by $C_n$.

\proclaim Lemma \at{lemt3}. Let $U$ be a $2$-dimensional vector space, and $C_n\simeq\P(U)\to\P(S^nU)$
be the $SL(U)$-equivariant embedding (whose   image is a rational normal curve).
 Let  $\C^*\subset SL(U)$ act  on  $\P(S^nU)$.  
We label the    $n+1$ fixed points $P_i, i=-n+2j, j = 0,\dots, n$ of the $\C^*$-action with an index proportional to its weight.    
 Then $P_{-n}$, $P_n$ lie 
on $C_n$ and $P_{-n+2j}=T^jP_{-n}\cap T^{n-j}P_n$, where $T^j$ denotes the 
$j$-dimensional osculating space to $C_n$.

\pr   We choose a coordinate system which diagonalizes the $\C^*$-action. 
Then the result follows by a straightforward computation.
\medskip

\proclaim Lemma\at{twodiff}. Let $S\in\S_{n,k}$ be a Steiner bundle. Let $Sym(S)^0$ be the 
connected component containing the identity of $Sym(S)$.
If there are two different one-parameter subgroups 
$\lambda_1,   
\lambda_2: \C^* \to   Sym(S)$ 
then $S$ is Schwarzenberger.  

\pr The proof is by induction on $k$.  If  $k=2$ the 
theorem is true by the theorem \rt{3.0}. By applying the theorem \rt{4.7}  to 
$\lambda_1$ we may suppose that 
the matrix representing $S$ is 
diagonal, and  that $H=\{x_0=0\}$ is  the fixed point  with minimum weight    of the 
dual action   $\lambda_1^*$ on 
$\P^{n\vee}$.
    By  \rt{elm} there is an elementary transformation
$$0\rig{}S'\rig{}S\rig{}\O_H\rig{}0$$
where  the matrix of $S'$
  is also diagonal (\re{elmm}), so that 
 $\lambda_1$ is a one-parameter subgroup of $Sym(S')$.  Let us suppose by 
contradiction that  $S$ is not Schwarzenberger; by the cor. \rt{cstar1}  
 we  find that $H$ is 
also the fixed point with minimum weight of the dual $\lambda_2^*$  (replacing  $\lambda_2$ with $\lambda_2^{-1}$ if necessary). Hence
by the same argument 
  also $\lambda_2$ is a one-parameter subgroup of $Sym(S')$, so that $S'$
is Schwarzenberger by the inductive assumption. It follows that $\lambda_1$ and 
$\lambda_2$ are contained in the same $SL(2) = Sym(S')$ and have the same two fixed points with 
minimum and maximum weight. By the lemma \rt{lemt3} $\lambda_1$ and $\lambda_2$ 
have the same fixed points and have also the same image in $SL(n+1)$. 
This is a contradiction.

{\it Proof of theorem \rt{t3}}

 Thanks to the theorem \rt{3.2}, the theorem \rt{t3} is equivalent to the following
 (the equivalence will be clear from the proof)

\proclaim Theorem\at{sym}. Let $S\in\S_{n,k}$ be a Steiner bundle. Let $Sym(S)^0$ be the 
connected component containing the identity of $Sym(S)$. Then
there is a $2$-dimensional vector space $U$ such that
$SL(U)$ acts over $\P^n=\P(S^nU)$ and according to this action
$Sym~(S)^0\subset SL(U)$. Moreover
$$Sym~(S)^0\simeq\left\{\matrix{0&\cr \C\cr \C^*\cr
SL(2)&\hbox{\ (this case occurs if and only if\ }S\hbox{\ is Schwarzenberger)
}}\right.$$

 We  prove the theorem 
\rt{sym}. The proof is by induction on $k$.  If  $k=2$ the 
theorem is true by the theorem \rt{3.0}. We may suppose that $G=Sym(S)^0$
has dimension $\ge 1$.
By the theorem \rt{t1} the matrix $A$ representing $S$ is triangulable.
By the prop. \rt{1.4} $W(S)$ is not empty and we pick up $H\in W(S)$.
By the corollary \rt{3.00} we may suppose that $W(S)$ is finite, hence $H$ is 
$G$-invariant. We repeat the argument of the proof of
the theorem \rt{3.2}. We get $\forall g\in G$ the diagram
$$\matrix{S&\rig{\phi}&\O_H\cr
\dow{i}\cr
g^*S&\rig{g^*\phi}&\O_H}$$
Since $h^0(S^*_{|H})\le 1$ we obtain that $\phi$ and $g^*\phi\circ i$ coincide
up to a scalar multiple. We get a commutative diagram
$$\matrix{0&\rig{}&S'&\rig{}&S&\rig{}&\O_H&\rig{}&0\cr
&&\dow{}&&\dow{\simeq}&&\dow{\simeq}\cr
0&\rig{}&g^*S'&\rig{}&g^*S&\rig{}&\O_H&\rig{}&0\cr}$$
It follows that $S'\simeq g^*S'$, hence $G\subset Sym(S')$ and by the 
inductive assumption $G\subset SL(U)$ and
$SL(U)$ acts over $\P^n=\P(S^nU)$. We remark that the above considered
elementary transformation
gives the decompositions $W=W'\oplus\C$, $I=I'\oplus\C$ 
such that the inclusion $Hom(W',V\otimes I')\subset
Hom(W,V\otimes I)$ identifies to the $SL(U)$-invariant inclusion
$S^{n+k-2}U\otimes S^nU\otimes S^{k-2}U\subset
S^{n+k-1}U\otimes S^nU\otimes S^{k-1}U$ according to the natural actions.
In fact no other morphism of $SL(U)$ in $SL(W)\times SL(S^nU)\times SL(I)$ can
give $S^{n+k-2}U\otimes S^nU\otimes S^{k-2}U$ as invariant summand of
$W\otimes S^nU\otimes I$.
Now consider the Levi decomposition $G=R\cdot M$ where $R$ is the
radical and $M$ is maximal semisimple. If $S$ is not Schwarzenbeger
we have $M=0$ and $G$ is solvable. By the Lie theorem $G$ is contained
(after a convenient basis has been fixed) in the subgroup 
$T=\left\{\left[\matrix{a&b\cr 0&{1\over a}
}\right] | a\in\C^*, b\in\C\right\}$. If there is a subgroup $\C^*$
properly contained in $T$ then there is a coniugate of $\C^*$ different from
itself and this is a contradiction by the lemma \rt{twodiff}.
If there is no subgroup $\C^*$ contained in $T$ then $G$ is isomorphic to
$\C\simeq\left\{\left[\matrix{1&b\cr 0&1
}\right] | b\in\C\right\}$.

\centerline{\bf References}

[AO] V. Ancona, G. Ottaviani, Stability of special instanton bundles on
$\P^{2n+1}$, Trans. AMS, 341, (1994) 677-693

[BS] G. Bohnhorst, H. Spindler, The stability of certain vector bundles on
$\P^n$, Lect. Notes Math. 1507, 39-50, Springer, 1992

[DK] I. Dolgachev, M.Kapranov, Arrangement of hyperplanes and vector bundles on
$\P^n$, Duke Math. J. 71 (1993), 633-664

[EH] Ph. Ellia, A. Hirschowitz, Voie Ouest I: g\'en\'eration de certains
fibr\'es sur les espaces projectifs et application, Journal Alg. Geometry
1 (1992)
531-547

[EP] D.Eisenbud, S.Popescu, The projective Geometry of the Gale transform,
alg-geom/9807127

[GKZ] I. M. Gelfand, M. M. Kapranov, A. V. Zelevinsky, Discriminants, resultants
and multidimensional determinants, Birkh\"auser, Boston 1994

[GKZ] I. M. Gelfand, M. M. Kapranov, A. V. Zelevinsky, Hyperdeterminants,
Adv. in
Math. 96 (1992), 226-263

[M] M. Maruyama, Elementary transformations in the theory of algebraic vector 
bundles, Lect. Notes Math. 961, 241-266, Springer 1982

[Schw] R.L.E. Schwarzenberger, Vector bundles on the projective plane,
Proc. London Math. Soc. 11 (1961), 623-640

[Simp] C. Simpson, Interactions between Vector bundles and Differential
equations,
problem list, preprint 1992, Europroj

[ST] H. Spindler, G. Trautmann, Special instanton bundles on $\P^{2n+1}$, their
geometry and their moduli, Math. Ann. 286 (1990), 559-592

[WZ] J. Weyman, A.V. Zelevinsky, Singularities of hyperdeterminants, Ann. Inst.
Fourier 46 (1996), 591-644

[V1] J. Vall\`es,   Fibr\'es de Schwarzenberger et coniques de droites
sauteuses,
preprint 1999

[V2]  J. Vall\`es,
Nombre maximal d'hyperplans instables pour un fibr\'e de Steiner, 
to appear in Math. Zeitschrift

\bigskip

Authors' addresses:

Dipartimento di Matematica

viale Morgagni 67/A

50134 FIRENZE ITALY

ancona@math.unifi.it

ottavian@math.unifi.it

 \end